\documentclass[12pt]{article}
\usepackage[utf8]{inputenc}
\usepackage[top=0.5in, left=0.65in, right=0.65in, bottom=0.75in]{geometry}
\usepackage{caption}
\usepackage{fancyhdr}
\usepackage{multicol}
\usepackage{graphicx}
\usepackage{dutchcal}
\usepackage[colorlinks=true,citecolor=blue,linkcolor=blue,urlcolor=blue]{hyperref}
\usepackage{amssymb, amsmath, amsfonts, amsthm, mathrsfs, tikz, tikz-cd, mathtools, mathdots, listings, etoolbox, enumerate, marvosym, wasysym}
\usepackage{verbatim}

\theoremstyle{definition}

\newcommand{\define}[4]{\expandafter#1\csname#3#4\endcsname{#2{#4}}}

\forcsvlist{\define{\DeclareMathOperator}{}{}}{im,coker,codim,rank,diam,vol,Tr,Nm,Ind,Res,Cor,Sym,Alt,Ob,Ext,Hom,Mor,End,Aut,Gal,lcm,gcf,sign,Spec,Proj,MaxSpec,Maxspec,Nil,JRad,dom,ext,Th,Bilin,tr,rad,eval,Frac,Supp,colim,res,id,Mat}

\forcsvlist{\define{\newcommand}{\mathsf}{}}{Sets,Grp,Ab,CRing,CAlg,Mod,Vect,Cat,Alg,Rep,Ring,Field,Rig,sym,fix}

\forcsvlist{\define{\newcommand}{\mathsf}{}}{WF,CSB,PP,WPP,ZF,CH,ZFA,ZFC,WPH,DC,CC,PA,PRA,NDS,RCA,ACA,WKL,ATR,LO,DLO,ACF,RCF,BT,BPIT,ZFCA,GCH}

\renewcommand{\AC}{\mathsf{AC}}

\forcsvlist{\define{\newcommand}{\mathbb}{}}{N,Z,Q,R,C,F}
\newcommand{\bR}{\mathbb{R}}

\newcommand{\bN}{\mathbb{N}}
\newcommand{\bQ}{\mathbb{Q}}
\newcommand{\bZ}{\mathbb{Z}}

\newcommand{\mc}[1]{\mathcal{#1}}

\newcommand{\cF}{\mc{F}}

\renewcommand{\geq}{\geqslant}

\renewcommand{\leq}{\leqslant}
\renewcommand{\epsilon}{\varepsilon}
\renewcommand{\emptyset}{\varnothing}
\newcommand{\isom}{\cong}

\newcommand{\set}[1]{\{#1\}}

\newcommand{\abs}[1]{\lvert#1\rvert}

\newcommand{\pset}{\mathcal{P}}

\newcommand{\es}{\emptyset}

\newcommand{\bbb}[1]{\ensuremath{\mathbb {#1}}}

\newcommand{\emp}{\varnothing}
\renewcommand{\phi}{\varphi}

\newcommand{\notarrow}{\kern .42em\not\kern -.42em\longrightarrow}

\usepackage{titlesec}
\titlespacing*{\section}{0pt}{*2.5}{*1}
\titlespacing*{\subsection}{0pt}{*2}{*1}

\begin{document}



    



























\begin{center}
    {\Large \textbf{A Gentle Introduction to the Axiom of Choice}} \vspace{2mm}
    
    {\large Andreas Blass and Dhruv Kulshreshtha}
\end{center}

\vspace{5mm}

\begin{multicols}{2}


\section{Introduction}

Described by David Hilbert as the axiom ``most attacked up to the present in the mathematical literature" \cite{Hilbert1926}, the axiom of choice ($\AC$) indeed has a fascinating history (see \cite{Moo82}). In 1883, Georg Cantor had proposed the well-ordering principle, i.e., that every set can be \textit{well-ordered}\footnote{We say that a set $X$ can be \textit{well-ordered} if there is a \textit{well-ordering} on $X$, i.e., a binary relation $\preceq$ on $X$ that is reflexive ($\forall x\in X, x\preceq x$), transitive ($\forall x,y,z\in X$, if $x\preceq y$ and $y\preceq z$, then $x\preceq z$), anti-symmetric ($\forall x,y\in X$, if $x\preceq y$ and $y\preceq x$, then $x=y$), total ($\forall x,y \in X$, $x\preceq y$ or $y\preceq x$), and well-founded (every nonempty subset of $X$ has a $\preceq$-minimal element). For example, the standard ordering $\leq$ on the set $\bN= \set{0,1,2,\ldots}$ of natural numbers is a well-ordering, and the standard ordering $\leq$ on the set $\bZ=\set{\ldots,-2,-1,0,1,2,\ldots}$ of integers fails only to satisfy the last condition, since, for example, $\bZ$ itself has no $\leq$-minimal element. (Note, however, that the integers can be well-ordered by, for example, constructing your favorite bijection of sets $h: \bN \to \bZ$ and ordering the integers by the induced ordering from $h$ and the standard ordering $\leq$ on $\bN$.) Similarly, the standard ordering $\leq$ on the set $[0,\infty)$ of non-negative real numbers---which itself has a $\leq$-minimal element---fails to satisfy the last condition since, for example, the subset $(0,\infty)$ of positive real numbers has no $\leq$-minimal element. (Unfortunately, one must work harder to prove that there is a well-ordering on $[0,\infty)$. Indeed, the initial motivation behind Cantor's well-ordering principle came from the question of whether the set $\bR$ of real numbers can be well-ordered---discussed in detail in \cite{Moo82}.)\label{footnote: well-ordered}}, without proof, as a valid law of thought. Since his so-called law was not immediately accepted, Cantor found himself seeking a proof. It was in order to provide such a proof that $\AC$ was first explicitly formulated by Ernst Zermelo in 1904 \cite{Zermelo1904}.

\medskip
\noindent\textbf{Axiom of Choice ($\AC$).} Given any family $\mc{F}$ of nonempty sets, there exists a function $f$ assigning to each member $A \in \cF$ an element $f(A)\in A$.\medskip


Such a function $f$, whose existence is given by $\AC$, is sometimes called a \textit{choice function} for $\mc{F}$, since $f$ can be thought of as ``choosing" an element from each set $A\in \mc{F}$. However, $\AC$ only guarantees the existence of such a function, without the means to construct one. 

More explicitly, $\AC$ does not say anything about one's ability to make or conceive these choices. Indeed, if $\cF$ is an infinite family, ``it is difficult to conceive how to make such choices---unless a rule is available to specify an element in each $A$" \cite{Moo82}. Moreover, in any instance of $\AC$ where there is available such a rule, $\AC$ is not actually needed.\footnote{More precisely, in any instance of $\AC$, where there is available such a rule that can be formalized in the first-order language of set theory, $\AC$ is not actually needed; other standard axioms---for example, Replacement and Pairing---will suffice.} A relevant example, due to Bertrand Russell, is discussed in Section \ref{sec: nonconstructivity}.

It is also worth explicitly noting that, even if an argument (in this case, an existence proof) seems to rely on making arbitrary ``choices," $\AC$ may not be needed. For instance, a standard proof of the Bolzano--Weierstrass theorem---every bounded sequence $(a_i)_{i\in \bN}$ in $\bR$ contains a convergent subsequence---relies on iteratively ``choosing" infinitely many terms $(a_{i_j})_{j\in J\subseteq \bN}$ from $(a_i)_{i\in \bN}$ satisfying some properties (see \cite[Theorem 2.5.5]{Abbott2015}). In this case, since the indexing set $\bN$ is well-ordered by $\leq$ (footnote \ref{footnote: well-ordered}), each ``choice" can be made definably by choosing the lowest-index term of the sequence satisfying those properties.

In this article, we begin with some common objections to $\AC$ (nonconstructivity and counterintuitive consequences). We then present three kinds of reasons to accept it. Although the present exposition is aimed at non-experts in set theory, we also include some lesser-known results. For more details about the axiom's history, mathematical significance, and consequences/equivalents see, for example, \cite{Jech73,Moo82,RubinRubin1985,HR98,Herrlich2006}. 

The footnotes in this article either define terminology or contain additional information that is not essential to the flow of the paper. A reader who is familiar with the relevant terminology can ignore them and still understand the paper.

We use the abbreviation $\ZFC$ for Zermelo--Fraenkel set theory with $\AC$---which is the standard axiomatic basis for mathematics---and the abbreviation $\ZF$ for Zermelo--Fraenkel set theory without assuming $\AC$. However, familiarity with the axioms of $\ZF$ is not assumed or necessary for this discussion.

\section{Objections}

In this section, we address two kinds of objections to $\AC$: its nonconstructive nature and (one of) its counterintuitive consequences. There are other objections---for example, the ability to make arbitrarily (hence uncountably) many choices---which are addressed in \cite{Moo82}.

\subsection{Nonconstructivity}\label{sec: nonconstructivity}

Nonconstructivity has been one of the primary objections to $\AC$.\footnote{Interestingly enough, many who raised this objection had themselves implicitly relied on $\AC$ and its consequences in their work \cite[Ch.1.7]{Moo82}.} To better explain the nonconstructivity of $\AC$, we describe a well-known example due to Bertrand Russell \cite{Russell1919}, replacing ``millionaire" from the original to account for inflation. As did Russell, we request ``a little goodwill" on the part of the reader in interpreting this example.

Imagine a billionaire who loves shopping for footwear, so he purchases pairs of shoes and socks until he possesses infinitely many pairs of each. One day, in a fit of eccentricity, he asks his butler to select and display one shoe from each pair. When his butler, accustomed to receiving precise instructions, asks how to decide which shoe to pick, the billionaire tells him to just choose the left shoe from each pair. The next day, in a similar fit, the billionaire asks his butler to select and display one sock from each pair. This time, when his butler asks how to decide which sock to pick, the billionaire is left dumbfounded: the socks in each pair are indistinguishable, so there is no analogous way to define a choice from each of the infinitely many pairs.

The point of Russell's anecdote is not to serve as a real life example. (Of course, in real life, one does not have infinitely many pairs of socks, so one's butler need not worry about having to choose one sock from each of infinitely many pairs. Furthermore, in real life, there may be physical factors distinguishing these socks, like slight differences in their weight.) Rather, the example is meant to highlight that it is in the latter of these scenarios that $\AC$ becomes relevant. Since there is no way to define a rule to choose one sock from each of the infinitely many pairs, one must rely on $\AC$ to guarantee the existence of a function that does the ``choosing."

\subsection{Mathematical inconveniences}\label{sec: BT and Div Paradox}
Often, a reason to not believe $\AC$ is its seemingly paradoxical consequences. If there were such a thing as a canonical\footnote{See, for example, ``the most spectacular [counterintuitive geometrical consequence of $\AC$]" in \cite[Ch.I]{Bell2009} and ``the most stunning [paradoxical result that is demonstrable in $\ZFC$ but not in $\ZF$]" in \cite[\S 5.2]{Herrlich2006}, as well as \cite[\S 9.4: Banach--Tarski paradox]{JustWeese1} and \cite[\S 1.3: A paradoxical decomposition of the sphere]{Jech73}.} example of a ``disastrous" consequence of the axiom of choice, it would be the Banach--Tarski paradox---which, roughly stated, says that a closed three-dimensional ball $B$ can be decomposed into finitely many pieces\footnote{In fact, five pieces suffice \cite{R-Robinson-BT-5-pieces-1947}; see \cite[Ch.\ 5]{Wagon1985} or \cite[Ch.\ 7]{Halbeisen-gentle-intro-2025}.}, which can be rearranged into two disjoint copies of $B$ using only rigid transformations.

More precisely, we say $B\approx C$ if there exist finite partitions $B = P_1\sqcup\cdots\sqcup P_n$ and $C=Q_1\sqcup\ldots\sqcup Q_n$ such that for each $1\leq i\leq n$, $P_i$ is congruent to $Q_i$.

\medskip\noindent\textbf{Banach--Tarski paradox ($\BT$).} Assuming $\AC$, given a closed three-dimensional ball $B$, there exists a decomposition $B=X\sqcup Y$ such that $X~\approx~B~\approx~Y$ \cite{BanachTarski1924}.\footnote{For a detailed discussion on $\BT$, see \cite{Wagon1985}.}\medskip

In spite of its seemingly counterintuitive nature, this so-called paradox is not a logical inconsistency, but rather a mathematical inconvenience.

$\BT$ (hence also $\AC$) implies another mathematically inconvenient result: the existence of subsets of $\bR^n$, for $n\geq3$, that are not \textit{(Lebesgue) measurable}.\footnote{For a formal and self-contained treatment of the Lebesgue measure on $\bR^n$, which is beyond the scope of this article, see \cite[Ch.1]{Stein-Shakarchi-3}. In the case of $\bR$, the intuition behind (Lebesgue) measurability comes from the ordinary notion of length: the measure of a point is zero and that of an open interval $(a,b)\subseteq \bR$ is its \textit{length}, denoted $m(a,b)=b-a$. Since any open subset $U$ of $\bR$ can be written uniquely as a disjoint union of open intervals \cite[Theorem 1.3]{Stein-Shakarchi-3}, the measure of $U$ can be defined as the sum of the lengths of these open intervals. Moreover, the Borel sets---defined as the smallest collection of subsets of $\bR$ that includes the open sets and is closed under taking complements and countable unions---are measurable \cite[Ch.1, \S3]{Stein-Shakarchi-3}. (For instance, since singletons are closed, hence Borel, all countable sets are Borel, hence measurable.) A subset $A$ of $\bR$ has \textit{measure zero} if it can be covered by arbitrarily small unions of open intervals: for any $\varepsilon>0$, there are open sets $U_i$ with $\sum_i m(U_i)<\varepsilon$ such that $A\subseteq \bigcup_i U_i$. (For instance, countably infinite sets have measure zero.) More generally, a set $A \subseteq \bR$ is \textit{measurable} if and only if it it differs from a Borel set by a set of measure zero \cite[Corollary 3.5]{Stein-Shakarchi-3}.\label{footnote: Lebesgue Measure}} More generally, as a consequence of $\BT$, there cannot exist any finitely additive measure that is defined on all subsets of $\bR^3$ and invariant under rigid transformations.

The standard example of a non-measurable subset of $\bR$ can be defined as follows. For $x,y\in \bR$, we write $x\sim_{\bQ} y$ to mean that $x-y \in \bQ$, where $\bQ$ denotes the set of rational numbers. This is an \textit{equivalence relation} on $\bR$.\footnote{A binary relation $\sim$ on a set $X$ is said to be an \textit{equivalence relation} if it is reflexive (footnote~\ref{footnote: well-ordered}), \textit{symmetric} ($\forall x,y\in X$, if $x\sim y$ then $y\sim x$), and transitive (footnote~\ref{footnote: well-ordered}). Given an equivalence relation $\sim$ on $X$, the $\sim$-\textit{equivalence class} of an element $x\in X$ is the set $[x]_\sim = \set{y\in X \mid x \sim y}$, and the \textit{partition} of $X$ by $\sim$ is the set $X/{\sim} = \set{[x]_\sim \mid x\in X}$ of all $\sim$-equivalence classes. Since each $\sim$-equivalence class is nonempty (as $x\in [x]_\sim$), $\AC$ guarantees the existence of a choice function on the partition $X/{\sim}$. We define a \textit{choice set} on $X/{\sim}$ as the range of such a choice function. Since distinct equivalence classes are disjoint, a choice set on $X/{\sim}$ will contain exactly one element from each equivalence class.\label{footnote: equiv rel and partiton}} A \textit{Vitali set} $V$ is a subset of $[0,1]$ containing exactly one element from each $\sim_\bQ$-equivalence class, or simply a choice set on $[0,1]/{\sim_\bQ}$ (see footnote~\ref{footnote: equiv rel and partiton}), whose existence is guaranteed by $\AC$. A standard argument of the non-measurability of such a $V$, relying on the translation-invariance of measure, can be found in \cite[Ch.1, \S3]{Stein-Shakarchi-3}.

Note, however, that the existence of a non-measurable set is not provable in $\ZF$, i.e., it is consistent with $\ZF$ (assuming the consistency of inaccessible cardinals) that all subsets of $\bR$ are measurable \cite{Solovay1970}.\footnote{For a brief discussion on assumptions that are consistent with $\ZF$ and incompatible with $\AC$, in the context of analysis, see \cite[\S14]{Schechter-handbook-analysis-1997}.} In this regard, it would seem as though the natural way of eliminating $\BT$ would be to insist that all subsets of $\bR$ (hence all subsets of $\bR^n$ \cite[Proposition 3.6]{Stein-Shakarchi-3}) are measurable and to suitably weaken $\AC$. However, this assumption leads to another paradox, which we describe briefly, that is just as unsettling as the one being eliminated.

For sets $A$ and $B$, we say $\abs{A}=\abs{B}$ if there is a bijection from $A$ onto $B$, and $\abs{A}\leq\abs{B}$ if there is an injection from $A$ into $B$. We write $\abs{A}<\abs{B}$ to mean $\abs{A}\leq\abs{B}$ and $\abs{A}\neq\abs{B}$.\footnote{For example, $\abs{\bN}=\abs{\bQ}$, $\abs{\bR}=\abs{(0,1)}$, $\abs{\bN}<\abs{\bR}$, and for any set $X$, $\abs{X}<\abs{\pset(X)}$. Here $\pset(X)$ denotes the \textit{power set} of $X$, i.e., the set of all subsets of $X$.} 

It is a consequence of $\AC$ that for any set $X$ and equivalence relation $\sim$ on $X$, $\abs{X/{\sim}} \leq \abs{X}$.\footnote{This statement, known as the \emph{partition principle} ($\PP$), can equivalently be stated as: ``if there is a surjection from $X$ onto $Y$, then there is an injection from $Y$ into $X$." If the injection were required to be a right inverse of the given surjection, this principle would be trivially equivalent to $\AC$. Without this additional requirement, however, equivalence between $\PP$ and $\AC$ remains one of the oldest and most frustrating open problems concerning $\AC$. See \cite{Moo82} for a detailed discussion.} This says that no partition of $X$ can strictly exceed $X$ (or even be incomparable with $X$) in size, which is in line with our intuition. With the additional assumption above, this falls apart:

\medskip\noindent\textbf{Division Paradox.} Assuming that all subsets of $\bR$ are measurable\footnote{The same conclusion can be derived from the assumption that all subsets of $\mathbb R$ have the \textit{Baire property}; see \cite[\S3]{TaylorWagon2019} and references there. This fact is reassuring if one doubts the consistency of inaccessible cardinals, because Shelah has shown that the consistency of ``all sets have the Baire property," unlike that of ``all sets are Lebesgue measurable," does not need the assumption that the existence of an inaccessible cardinal is consistent \cite{Shelah1984-inaccessible}.}, we have $\abs{\bR/{\sim_\bQ}}>\abs{\bR}$ (due to Sierpi\'nski and Mycielski; see \cite[\S3]{TaylorWagon2019} and references there).\medskip

That is, $\bR$ is partitioned into strictly more parts, under the $\sim_\bQ$ equivalence relation, than there are real numbers! This equally counterintuitive result, among others, is why ``we have learned to live with nonmeasurable sets" \cite{TaylorWagon2019}.

What is remarkable is that the same partition $\bR/{\sim_\bQ}$ of the reals that gives a non-measurable set in the presence of $\AC$ is also responsible for the division paradox when all sets of reals are assumed to be measurable. 

\section{Why use some Choice?}
In this section, we give three (kinds of) arguments in favor of using $\AC$, or at least some of its consequences.

\subsection{Because a lot of math follows from $\AC$ or needs $\AC$}\label{sec: math follows from AC}

Thus far, we have seen some uses of $\AC$ that are not necessarily intuitive. For instance, in the introduction, we mentioned that $\AC$ was used to prove the well-ordering principle. In fact, this use is necessary, in the sense that $\AC$ is equivalent to the well-ordering principle. $\AC$ is also regularly used throughout mathematics in another equivalent form, namely Zorn's Lemma\footnote{``Zorn's Lemma" is the standard name for this principle. A more accurate title would include the names of Kuratowski and Hausdorff, among others. See \cite{Campbell-origin-of-Zorn-1978} for details on the history and origins of the various maximal principles equivalent to $\AC$.}: if $(P,\preceq)$ is a nonempty \textit{partially ordered} set in which every nonempty \textit{chain} $C\subseteq P$ has a $\preceq$-upper bound, then $P$ has a $\preceq$-maximal element.\footnote{A binary relation $\preceq$ on $P$ is said to be a \textit{partial order} if it is reflexive, anti-symmetric, and transitive (see footnote~\ref{footnote: well-ordered}). A subset $C$ of a partially ordered set $(P,\preceq)$ is said to be a \textit{chain} if the restriction of $\preceq$ to $C$ is total (again, see footnote~\ref{footnote: well-ordered}).} The equivalence of these statements may not seem obvious at all, as reflected by a famous (ironic) quote of Jerry Bona from 1977: ``The Axiom of Choice is obviously true; the Well-Ordering principle is obviously false; and who can tell about Zorn's Lemma?" \cite{Schechter-handbook-analysis-1997}. The curious reader, who at this point has surely attempted showing that these statements are equivalent, may find a proof in \cite[Theorem 2.1]{Jech73}.

In this subsection, we provide results that rely on $\AC$ or its consequences and are more intuitive (or at least closer to the interests of a typical non-logician) than those above. We draw from introductory analysis and algebra, and refer the curious reader to more detailed sources for similar discussions in other fields.

In introductory real analysis, one comes across two notions of a function $g:\bR\to\bR$ being ``continuous" at a point $x\in \bR$. $g$ is said to be \textit{continuous} at $x$, if for any $\varepsilon>0$, there exists $\delta>0$, such that for any $y\in \bR$, if $\abs{x-y}<\delta$ then $\abs{g(x)-g(y)}<\epsilon$; and $g$ is said to be \textit{sequentially continuous} at $x$ if for any sequence $(a_n)_{n\in\bN}$ of real numbers converging to $x$, the sequence $(g(a_n))_{n\in\bN}$ converges to $g(x)$. 

Assuming $\AC$, these notions agree: $g$ is continuous at $x$ if and only if $g$ is sequentially continuous at $x$. However, the backward direction of this equivalence cannot be proved in $\ZF$, i.e., it is consistent with $\ZF$ that the backwards direction fails (\cite{Jaegermann1965}; discussed in \cite[Ch.1.2]{Moo82}).

A natural question to raise is whether we need the ``full strength" of $\AC$ for this equivalence to hold. This turns out not to be the case: the backwards direction is equivalent to countable choice from sets of reals, i.e., the assertion that any countable family of nonempty subsets of $\bR$ admits a choice function, abbreviated $\CC(\bR)$ \cite[Theorem 4.54, 8 and 10]{Herrlich2006}. $\CC(\bR)$ is the restriction of $\AC$ where $\cF$ is countable, i.e., $\abs{\cF}\leq\abs{\bN}$, and each element of $\cF$ is a nonempty subset of $\bR$. This is an obvious consequence of $\AC$ restricted to countable families of arbitrary nonempty sets, abbreviated $\CC$, which itself is strictly weaker than $\AC$.\footnote{Both of the implications $\AC\implies \CC\implies \CC(\bR)$ are strict. See \cite[Ch.8]{Jech73} for the irreversibility of the first implication. Since $\CC(\bR)$ holds in every \textit{permutation model} and is a \textit{boundable statement} (see \cite[III.\S2]{HR98} and \cite[Ch.6]{Jech73})---hence, transferrable to a $\ZF$ model---the irreversibility of the second implication can be seen through any permutation model where $\CC$ fails, such as the Basic Fraenkel Model (see \cite[III.\S2.1]{HR98} or \cite[Ch.4.3]{Jech73}).}

Note, furthermore, that $\CC(\bR)$ implies that the Lebesgue measure (see footnote~\ref{footnote: Lebesgue Measure}) is \textit{countably additive}, i.e., if $(A_i)_{i\in\bN}$ are pairwise disjoint measurable subsets of $\bR$, then $m(\bigcup_iA_i) = \sum_im(A_i)$ \cite[\S5.1(E13)]{Herrlich2006}. Countable additivity is not provable in $\ZF$, i.e., it is consistent with $\ZF$ that the Lebesgue measure is not countably additive.\footnote{For example, in the Feferman--L\'evy Model (\cite{Feferman-Levy-1963}; see \cite[III.\S1.9]{HR98}), $\bR$ is a countable union of countable sets, so there is no nontrivial countably additive measure on $\bR$. (Of course, this in turn means that $\CC(\bR)$ fails in this model, and so does the aforementioned equivalence of continuity notions.)}

In this regard, it may seem as though some countable analog of $\AC$ is sufficient for many purposes in real analysis.\footnote{Indeed, some analysts who were critical of $\AC$---while relying on its consequences implicitly in their own work---were willing to accept its restriction to countable families (see \cite[\S1.7 and Appendix 1]{Moo82}).} This is far from the case in algebra: ``Algebraists insisted that [$\AC$], whether in the guise of Zorn's Lemma or the Well-Ordering Theorem, had become indispensable to their discipline" \cite{Moo82}. We mention a few instances of this.

Let $F$ be a field and $V$ be an $F$-vector space. If $V$ has a finite spanning set, i.e., if there are $x_1,\ldots,x_k\in V$ such that any $y\in V$ can be written as a linear combination $y=\sum_{i=1}^k c_ix_i$, with $c_i\in F$, then it is not hard to show, without appealing to $\AC$, that $V$ has a (finite) basis. Furthermore, if an $F$-vector space $V$ has a finite basis, it is once again not hard to show, without appealing to $\AC$, that any two bases of $V$ must have the same cardinality---so the \textit{dimension} of $V$ can be defined as the cardinality of any basis of $V$.

On the other hand, there are plenty of vector spaces that do not have finite spanning sets.\footnote{For example, the space $\bR[x]$ of polynomials in $x$ with real coefficients, the space of continuous functions from $\bR$ to $\bR$, the space of sequences $a_n:\bN\to \bR$, all considered as $\bR$-vector spaces, as well as the set $\bR$ considered as a $\bQ$-vector space.} In the general (not necessarily finite) case, it turns out that $\ZF$ itself is not strong enough to guarantee the existence of a basis or that bases have the same cardinality. That is, it is consistent with $\ZF$ that there is a vector space with no basis, and that there is a vector space with two bases of different cardinalities (\cite{Lauchli1962}; see \cite[Theorem 10.11 and Ch.10.3(5)]{Jech73}). This means that in $\ZF$, the notion of dimension may not be well-defined. Furthermore, in contrast to the continuity notions, we can no longer get away with weaker forms of $\AC$: it is equivalent to $\AC$ that every vector space has a basis \cite{Blass1984}.

The full strength of $\AC$ is also needed in other areas of algebra. For instance, it is equivalent to $\AC$ that every (nonzero) commutative ring with unit has a maximal (proper) ideal \cite{Krull1929,Hodges1979}; and that for every abelian group $G$ and subgroup $H$ of $G$, there is a set of representatives for the quotient $G/H$ \cite{Keremedis-1998}.

For results and discussions on results in other fields of mathematics that rely on $\AC$ or its consequences, we refer the interested reader to \cite[Ch.10]{Jech73}, \cite{Hodges1974}, \cite{Moo82}, \cite{Galvin-Komjath-1991}, \cite{HR98} and references there, \cite[\S4]{Herrlich2006}, and \cite{Halbeisen-et-al-AC-Rings-2019}.

\subsection{Cracking walnuts with a sledgehammer}

In this subsection, we briefly describe how sometimes cracking walnuts with a sledgehammer can result in simpler proofs.\footnote{The ``cracking walnuts with a hammer" analogy is due originally to Alexander Grothendieck \cite{Grothendieck-1985-1987}; see \cite{Mclarty-Grothendieck-2007} for a brief description. The phrase ``cracking walnuts with a sledgehammer" was once used to describe the second author's (unintentional) use of unnecessarily powerful mathematical machinery in place of the much simpler intended proof, while the author was a course assistant for an introductory mathematical writing course. The first author has similarly been accused of ``using a bazooka to kill an ant."} More formally, in this subsection, we discuss some results that are provable in $\ZF$ but have simpler proofs using $\AC$.

We discuss a few results about cardinalities of sets and a classical result in combinatorics, providing references to standard $\ZF$ proofs where required.

Recall that for sets $A$ and $B$, we say $\abs{A}=\abs{B}$ if there is a bijection from $A$ onto $B$, and $\abs{A}\leq\abs{B}$ if there is an injection from $A$ into $B$. Additionally, for any set $X$ and $n\in \bN_{>0}$, we write $n\times X$ to denote the cartesian product $\set{0,1,\ldots,n-1}\times X$.

\medskip\noindent\textbf{Division by $m$.} For any sets $A$ and $B$, and any $m\in\bN_{>0}$, if $\abs{m\times A}\leq \abs{m\times B}$, then $\abs{A}\leq\abs{B}$ (Lindenbaum [unpublished] and Tarski \cite{Tarski-1949}).\footnote{This result has its own interesting history, which is summarized in \cite{DC94}.}
\medskip

For finite sets $A$ and $B$, this follows from straightforward arithmetic.
It is when dealing with infinite sets that this becomes a fairly nontrivial argument over $\ZF$ (see \cite{DC94}).

On the other hand, cardinal arithmetic is very tame in the presence of $\AC$ \cite[Ch.5]{Jech-set-theory-1997}. Indeed, assuming $\AC$, for any nonempty sets $X$ and $Y$, if at least one is infinite, then $\abs{X\times Y}=\max\set{\abs{X},\abs{Y}}$. 
This easily tackles the infinite case of the preceding theorem. However, this specific argument heavily relies on $\AC$: even if we just assume that for every infinite set $X$, $\abs{X\times X}=\abs{X}$, then $\AC$ follows (\cite{Tarski-1924a}; see \cite[Ch.11]{Jech73}). Moreover, $\AC$ is also equivalent to the assertion that for any sets $X$ and $Y$, $\abs{X}\leq\abs{Y}$ or $\abs{Y}\leq\abs{X}$ \cite{Hartogs1915}, which is necessary in defining $\max\set{\abs{X},\abs{Y}}$.
\medskip

For the next result, we introduce the following notation: given a set $X$, let $\mathcal{W}(X)$ be the set of all well-orderable subsets of $X$, i.e., those subsets of $X$ that can be well-ordered (see footnote \ref{footnote: well-ordered}).

\medskip\noindent\textbf{More well-orderable subsets.} For any set $X$, $\abs{X}<\abs{\mathcal{W}(X)}$ (Tarski \cite{Tarski-1939}).
\medskip

Recall that $\AC$ is equivalent to the well-ordering principle. So, over $\ZF$---where some sets may admit no well-ordering---this is a strengthening of Cantor's theorem, which states that $\abs{X}<\abs{\mc{P}(X)}$, where $\mc{P}(X)$ denotes the power set of $X$, i.e., the set of all subsets of $X$. By contrast, assuming $\AC$, $\mathcal{W}(X)=\mc{P}(X)$, so this result follows from Cantor's theorem.
\medskip

Another classical $\ZF$ result on cardinalities is:

\medskip\noindent\textbf{Cantor--Schr\"oder--Bernstein~Theorem~($\CSB$).} For any sets $X$ and $Y$, if $\abs{X}\leq \abs{Y}$ and $\abs{Y}\leq \abs{X}$, then $\abs{X}=\abs{Y}$.\footnote{This result also has its own extremely interesting history, which is discussed in detail in \cite{Hinkis-CSB-book-2013}. In particular, a more accurate title for the theorem would also include the names of Dedekind and Zermelo.}
\medskip

Over $\ZF$, the proof is fairly hands-on. Given injections $f:X\to Y$ and $g:Y\to X$, rather than merely proving the \textit{existence} of a bijection, one must \textit{construct} a bijection $h:X\to Y$ by stringing together these maps (see, for example, \cite{Jech-set-theory-1997}). Assuming $\AC$, the sets $X$ and $Y$ are well-orderable. In this case, standard results about well-orderings provide a ``conceptually simpler" proof.\footnote{In this case, the ``simplicity" of the $\AC$-proof of $\CSB$ depends on taking for granted standard facts about well-orderings, so it may not necessarily be mathematically simpler. For instance, this argument relies on the trichotomy of well-orderings, i.e., given two well-orderings $W_1$ and $W_2$, either  $W_1\isom W_2$ or one is isomorphic to an initial segment of the other, and on the fact that if $W$ is a well-ordered set and $V$ is a subset of $W$ with the induced ordering, then $W$ is not isomorphic to a proper initial segment of $V$. In contrast, the $\ZF$ proof, though tedious, does not require this additional machinery. A discussion by Hamkins and Schweber on whether or not this proof is actually simpler is in the comments of Hamkins' answer in \cite{Hamkins-Schweber-MathOverflow-simpler-proofs-AC}.}  
\medskip

Our final example is Hindman's theorem from combinatorics, which states that for any finite coloring of $\bN$, there is an infinite $X\subseteq \bN$, all of whose finite nonempty sums of distinct elements are assigned the same color. To make this precise, for $A\subseteq \bN$, we define $\mathrm{FS}(A) = \set{\sum_{a\in F} a \mid F\neq \es \text{ is a finite subset of } A}$. The theorem can now be stated as follows.

\medskip\noindent\textbf{Hindman's Theorem.} For any finite partition $\bN = C_1\sqcup\dots\sqcup C_n$, there exists an infinite subset $X\subseteq \bN$ and $i\leq n$ such that $\mathrm{FS}(X) \subseteq C_i$.
\medskip

Hindman's original proof \cite{Hindman-Theorem-original-1974} is purely combinatorial and uses fairly elementary machinery: it can be formalized in second-order arithmetic with room to spare, as shown by Hirst and published in \cite{BlassHirstSimpson1987}. In particular, no use is made of $\AC$. On the other hand, the proof is quite difficult to follow. Indeed, Hindman has himself suggested that the original proof can be used as a torture device for graduate students (see \cite{Goldbring2022}).

A standard alternative is the Galvin--Glazer argument, described by Hindman as ``the shortest and, to my mind, prettiest proof\," of Hindman's theorem \cite{Hindman-Theorem-expository-1979}. Formalizing this proof requires several more levels of the cumulative hierarchy, and this proof also relies on two applications of $\AC$. However, in contrast to the original, this argument is one that can reasonably be remembered. Since this proof is still outside of the scope of this article, we provide a brief outline of the applications of $\AC$ in this proof, while omitting details. The detailed proof can be found in \cite{Hindman-Theorem-expository-1979}.

The set of \textit{nonprincipal ultrafilters} on $\bN$ has a natural topology, which makes it a compact Hausdorff space. One can define a \textit{semigroup} structure on this space (see \cite{Hindman-Theorem-expository-1979} or \cite{Goldbring2022}). The first application of $\AC$ comes in the form of a special case of the Boolean Prime Ideal Theorem (Form 14 in \cite{HR98}) to prove the existence of a nonprincipal ultrafilter on $\bN$. A second application of $\AC$ is in the form of Zorn's lemma to get a minimal nonempty closed subsemigroup---which provably contains just one ultrafilter that is \textit{idempotent} (with respect to the semigroup operation). A combinatorial argument then gives a set $X$, as required in the theorem, without any additional use of $\AC$.

\subsection{Because you can!}
Is $\AC$ safe to use? Or might $\ZFC$ be inconsistent? Equivalently, might the negation of $\AC$ (abbreviated $\neg\AC$) be a theorem of $\ZF$? In view of the doubts about $\AC$ expressed by prominent mathematicians \cite{Moo82}, such questions should be taken seriously.

Another (opposite) issue is whether $\AC$ is redundant. Might it be a theorem of $\ZF$? Equivalently, might the theory $\ZF + \neg\AC$ be inconsistent?

The answers to both questions are on the side of consistency. That is, each of $\ZFC$ and $\ZF + \neg\AC$ is consistent. But to prove these answers we need an additional assumption, namely that $\ZF$ (without any commitment for or against $\AC$) is consistent. The additional assumption is
required because of G\"odel's second incompleteness theorem \cite{Godel1931}, which says that a reasonable\footnote{``Reasonable'' means that (1)~the theory is consistent, (2)~it can prove basic facts about the arithmetic of natural numbers, and (3)~its axioms can be listed by an algorithm. In particular, $\ZFC$ is reasonable.} theory cannot prove its own consistency. Thus, what is actually proved is that if $\ZF$ is consistent then so are $\ZFC$ and $\ZF + \neg\AC$. In other words, if $\ZF$ were to become inconsistent when one adds $\AC$ (or $\neg\AC$) as a new axiom, then $\ZF$ itself was already inconsistent; the inconsistency is not the fault of the added axiom.

The goal of this subsection is to describe the ideas that underlie the proofs of these consistency statements. That will require some preliminary work, describing the main intuition that underlies $\ZF$, namely the \textit{universe} $V$ obtained as a cumulative hierarchy of sets.

This hierarchy is built as follows. Begin with a collection of objects that are not sets; these are called \emph{atoms} or \emph{urelements}. They constitute level~0 of the hierarchy. Level~1 consists of all sets of atoms. Level~2 consists of all the new sets that can be formed using atoms and sets of atoms as members. Continuing in this way, each level consists of all the new sets whose elements are at strictly earlier levels. This process is to be continued not only through level~$n$ for all natural numbers $n$, but transfinitely forever.

This description of the cumulative hierarchy involves vague notions, ``all sets'' and ``transfinitely forever,'' so it cannot be used directly as a basis for mathematical proofs. Despite the vagueness, the intuition of the cumulative hierarchy supports some precise statements about sets, and a useful collection of those was isolated in an axiomatic system $\ZFA$ (meaning $\ZF$ with atoms, and $\ZFCA$ adds $\AC$). $\ZF$ is equivalent to $\ZFA$ plus the axiom that there are no atoms\footnote{Even without atoms, there are plenty of sets: $\emp$ at level~1, $\{\emp\}$ at level~2, $\{\{\emp\}\}$ and $\{\emp,\{\emp\}\}$ at level~3, and so on. Infinite sets appear at the transfinite levels.}, so $\ZFA$ is consistent if $\ZF$ is. More importantly for our purposes, the consistency of $\ZF$ implies the consistency of $\ZFA$ plus the axiom that there are infinitely many atoms \cite[Problem 4.6.1.]{Jech73}. It turns out that atoms are never really needed in mathematics, so $\ZFC$, rather than $\ZFCA$, has become the widely accepted axiomatic basis for set theory and thus for the rest of mathematics.

The proofs of consistency of $\AC$ and of its negation were developed in three phases, to be described below.  In all phases, one starts with a universe $V$, satisfying all axioms of $\ZF$ (or $\ZFA$ if there are atoms)\footnote{We think of $V$ as consisting of all the atoms and all the sets in the cumulative hierarchy. Note that this $V$ is not a set; such collections are often called proper classes. For technical reasons, many authors prefer to start with a set-sized model of $\ZF$, and they often impose additional conditions on it. Those technicalities will not affect our discussion, so we omit them and work with the ``whole'' universe $V$.}, and one modifies it carefully to satisfy $\AC$ or to satisfy $\neg\AC$, while still satisfying the axioms of $\ZF(\mathsf{A})$. The relevant modifications are quite different in the three phases.

The first phase began in 1922 with Abraham Fraenkel's proof that $\AC$ cannot be proved in $\ZFA$ \cite{Fraenkel-1922-indep-AC}. He began with a universe satisfying $\ZFA$ and having infinitely many atoms, and he defined a sub-universe $S$ of ``sufficiently symmetric"\footnote{Here, ``symmetry" refers to behavior when the atoms are arbitrarily permuted. Such a permutation $\pi$ produces a permutation of $V$. (In fact, $\pi$ produces an automorphism of $V$.) To apply $\pi$ to a set, simply apply $\pi$ to all its elements. Since the elements occur earlier in the cumulative hierarchy, this is a legitimate inductive definition. A set or atom $x$ is ``sufficiently symmetric" if there is a finite set $E\subseteq A$ such that all permutations that fix all elements of $E$ also fix $x$. (It is tempting to simplify this definition by requiring $x$ to be fixed by all permutations. Unfortunately these ``completely symmetric" sets don't satisfy $\ZFA$. For example, both $\emp$ and $A$ are completely symmetric but they are unequal despite having the same completely symmetric members, namely none.)} sets to violate $\AC$---for example, the set of atoms admits no symmetric well-ordering (or even linear ordering)---while satisfying all of the $\ZFA$ axioms. Fraenkel's sub-universe consists of the hereditarily sufficiently symmetric sets and atoms, i.e., those $x$ such that $x$ and all its members, and all their members, $\dots$ are sufficiently symmetric. It is fairly easy  to see that $\AC$ is false in this sub-universe. For example, all the atoms and the set $A$ of atoms are sufficiently symmetric, but no linear ordering of $A$ (let alone a well-ordering) is sufficiently symmetric. More work is needed to show that this sub-universe satisfies all the $\ZFA$ axioms; see \cite{Jech73} or \cite{Halbeisen-gentle-intro-2025}.

Fraenkel's method was extended, especially by Andrzej Mostowski, to obtain considerable information about non-implications between various consequences of $\AC$ \cite{Fraenkel-1937,Mostowski-1939}. For example, even if one adds to $\ZFA$ an axiom saying that every set can be linearly ordered, one still cannot deduce the axiom of choice.  Later, Ernst Specker systematized the theory of such \emph{permutation models} (or Fraenkel--Mostowski models or Fraenkel--Mostowski--Specker models) in terms of an arbitrary group of permutations of atoms (to define ``symmetric'') and an arbitrary normal filter of subgroups (to define ``sufficiently'') \cite{Specker-1957}. For more information about permutation models of $\ZFA$, see \cite[Chapter~4]{Jech73}.

Unfortunately, the method of permutation models depends crucially on the presence of infinitely many atoms. It tells us nothing about \textit{pure sets}\footnote{Pure sets are sets that have no atoms among their members, members of members, $\dots$.}, which are completely symmetric, hence survive the shrinking from $V$ to $S$. Since, for example, the usual constructions of the real numbers make them pure sets, the question whether \bbb R can be well-ordered has the same answer in $S$ as in $V$.  Even after the development of permutation models, one could reasonably imagine that $\AC$ is provable from $\ZF$.  Eliminating that possibility would have to wait 40 years for the third phase.

The second phase is mostly work of Kurt G\"odel, who proved in \cite{Godel1938} that, if $\ZF$ is consistent, then so is $\ZFC$. He showed how to find, within any universe satisfying $\ZF$, a sub-universe $L$, called the \emph{constructible universe}, satisfying $\ZFC$. The definition of $L$ proceeds similarly to the cumulative hierarchy that produces the universe $V$ of all sets. The only difference is that, at each level, one puts into that level only the \emph{definable} subsets
of the union of the previous levels.\footnote{Here ``definable'' means first-order definable with parameters from earlier levels. (This notion of definability was less known at the time, 1938, and in his book \cite{Godel-1940} G\"odel circumvented it by a set of operations on sets, now called the G\"odel operations, that parallel the construction of first-order formulas.)}

Recall that $\AC$ is used to produce sets that cannot be defined; definable sets are provided by the other axioms. So it should seem strange that $\AC$ is true in a universe $L$ produced entirely from definable sets.\footnote{It is also worth noting that not everything in $L$ is entirely definable. The ordinal numbers that index the levels of the constructible hierarchy are the same as in the cumulative hierarchy; they are not subject to any definability restrictions.} This situation arises from using definitions themselves in order to make choices. One well-orders sets in $L$ according  to (1)~the level at which they enter the constructible hierarchy, (2)~in case of a tie, the lexicographic order of their defining formulas, and (3) in case it's still a tie, the relative order of the parameters in the definition, as given by well-ordering the previous levels.

In the same work, G\"odel proved that $L$ satisfies, in addition to $\ZFC$, the generalized continuum hypothesis ($\GCH$), i.e., the statement that, for any infinite cardinal $\kappa$, the cardinality $2^\kappa$ of its power set is the first cardinal greater than $\kappa$. In particular, the cardinality of \bbb R is the first uncountable cardinal $\aleph_1$.\footnote{Subsequent work by Jensen \cite{Jensen-1972} showed that numerous other useful combinatorial principles hold in $L$. Jensen himself showed that Souslin's Hypothesis (i.e., an affirmative answer to Souslin's question \cite{Souslin-1920}) is not provable in $\ZFC$ (not even in $\ZFC+\GCH$) because it is false in $L$. Jensen's combinatorial principles have found applications in topology, algebra, and functional analysis.}

The third phase of consistency results about $\AC$ began with Paul Cohen's invention of the method of forcing \cite{Cohen1963,Cohen1964}.\footnote{For more modern approaches, see \cite{Shoenfield-1971}, \cite{Jech73}, \cite{Jech-set-theory-1997}, or \cite{Kunen1980}.} The key idea here is, starting with $V$ satisfying $\ZFC$, to produce a larger model, having $V$ as  a substructure and satisfying $\ZF$ and the negation of $\AC$.

But we were working with a universe $V$ containing all sets; how can one add more? Forcing enables one to consider a generalized notion of set. For these generalized sets, statements like $x\in y$ are not simply true or false; they can have other truth values\footnote{The following description is not Cohen's original version but an equivalent simplification invented by Dana Scott and Robert Solovay [unpublished]. Standard references for it include \cite{Jech-set-theory-1997,Bell-2011}.}. The truth values are the elements of some complete Boolean algebra $B$. Once $B$ is chosen, one can build the generalized cumulative hierarchy, adding at each level all generalized sets whose generalized elements are at earlier levels. Cohen proved that the resulting generalized universe $V^B$ satisfies $\ZFC$; its other properties depend on the choice of $B$. Note that the proof that $\ZFC$ is satisfied, i.e., all its axioms have the truth value $1$ (the top element of $B$), is complicated, because the very notion of set has been changed by allowing truth values from $B$.

This construction enabled Cohen to prove the independence of $\GCH$ from $\ZFC$; that is, he designed $V^B$ to violate $\GCH$. Specifically, it had $2^{\aleph_0}=\aleph_2$. Subsequent work has used this method to prove consistency of a  great variety of propositions across many fields of mathematics.

But we wanted a model of $\ZF$ that violates $\AC$, and $V^B$ doesn't accomplish that. If $V$ satisfies $\ZFC$ then so does $V^B$. For violations of $\AC$, one needs to import the idea of symmetry from Fraenkel's work in the first phase. Now, symmetry can no longer be based on permutations of atoms, since we have no atoms. Instead, symmetry is based on automorphisms of the Boolean algebra $B$. Any automorphism $\pi$ of $B$ produces, by induction on levels, an ``automorphism'' $\bar\pi$ of $V^B$. The quote marks are because we can't just apply $\bar\pi$ to the generalized sets in $V^B$; we must simultaneously apply $\pi$ to the Boolean values in~$B$. Once this is taken into account, we can use an arbitrary group of automorphisms of $B$ and an arbitrary normal filter of subgroups, just as in the first phase, to select a sub-universe of $V^B$, called a \emph{symmetric model}. An important theorem is that such models satisfy $\ZF$.

For suitable choices of $B$, the group, and the filter, the symmetric model violates $\AC$. The choices here influence which consequences of $\AC$ hold in the symmetric model and which fail.  For example, in some symmetric models all sets can be linearly ordered; in others they cannot. In applications of Cohen's method, a crucial part of the work is finding a suitable $B$, and this can often be simplified. Any complete Boolean algebra $B$ is completely specified when a dense subset $D$ is
given. Here ``dense'' means that every non-zero element of $B$ is above a  non-zero element of $D$, and ``given'' means that we have the set $D$ and its partial ordering induced from $B$. It is very often easier to find and describe a suitable partial order $D$  than to describe the corresponding complete Boolean algebra. Furthermore, the theory of Boolean-valued models can be rewritten in terms of $D$, without even mentioning $B$. This way of using a dense set $D$ leads back to (a variant of) Cohen's original forcing method; his ``conditions'' constitute such a dense set.

\bibliography{References}
\bibliographystyle{alpha}

\end{multicols}

\end{document}